\documentclass[11pt]{article}
\usepackage{psfig}
\usepackage{amssymb}
\newcommand{\be}{\begin{equation}}
\newcommand{\ee}{\end{equation}}
\newcommand{\bea}{\begin{eqnarray}}
\newcommand{\eea}{\end{eqnarray}}
\newcommand{\binom}[2]{{#1 \choose #2}}

\begin{document}

%\DeclareGraphicsExtensions{.jpg,.pdf,.png}

\title{Restricted Partition Functions as \\
Bernoulli and Euler Polynomials of Higher Order}
\author{Boris Y. Rubinstein${}^{\dag}$ and Leonid G. Fel${}^{\ddag}$\\
\\
${}^{\dag}$Department of Mathematics, University of California, Davis,
\\One Shields Dr., Davis, CA 95616, U.S.A. \\
and \\
${}^{\ddag}$Department of Civil and Environmental Engineering,\\
 Technion, Haifa 32000, Israel}
\date{\today}

\maketitle
%%%%%%%%%%%%%%%%%%%%%%%%%%%%%%%%%%%%%%%%%%%%%%%%%%%%%%%%%%%%%%%%%%%%%
\begin{abstract}
Explicit expressions for restricted partition
function $W(s,{\bf d}^m)$ and its quasiperiodic components
$W_j(s,{\bf d}^m)$ (called {\em Sylvester waves}) for a set of positive
integers ${\bf d}^m = \{d_1, d_2, \ldots, d_m\}$ are derived.
The formulas are represented in a form of a finite sum over
Bernoulli and Euler polynomials of higher order with periodic coefficients.
A novel recursive relation for the Sylvester waves is established.
Application to counting algebraically independent homogeneous
polynomial invariants of the finite groups is discussed.
% polynomial invariants of degree $s$ of the finite group $G$ is discussed.
\end{abstract}
%%%%%%%%%%%%%%%%%%%%%%%%%%%%%%%%%%%%%%%%%%%%%%%%%%%%%%%%%%%%%%%%%%%%%
%\newpage
\section{Introduction}
\label{intro}
The problem of partitions of positive integers has long history started from
the work of Euler who laid a foundation of the theory of partitions 
\cite{GAndrews}, introducing the idea of generating functions.
Many prominent mathematicians contributed to the development of the theory 
using the Euler idea.

J.J. Sylvester provided a new insight and made a
remarkable progress in this field. He found \cite{Sylv1,Sylv2} the 
procedure enabling to determine a {\it restricted} partition functions, 
and described symmetry properties of such functions. The restricted
partition function $W(s,{\bf d}^m) \equiv W(s,\{d_1,d_2,\ldots,d_m\})$ is a
number of partitions of $s$ into positive integers  $\{d_1,d_2,\ldots,d_m\}$,
each not greater than $s$. The generating function for $W(s,{\bf d}^m)$ 
has a form
\be
F(t,{\bf d}^m)=\prod_{i=1}^m\frac{1}{1-t^{d_{i}}}
 =\sum_{s=0}^{\infty} W(s,{\bf d}^m)\;t^s\;,
\label{genfunc}
\ee
where $W(s,{\bf d}^m)$ satisfies the basic recursive relation
\be
W(s,{\bf d}^m) - W(s-d_m,{\bf d}^m) = W(s,{\bf d}^{m-1})\;.
\label{SW_recursion}
\ee
Sylvester also proved the statement about splitting of the partition 
function into periodic and non-periodic parts and showed that the 
restricted partition function may be presented as a sum of "waves", which 
we call the {\em Sylvester waves}
\be
W(s,{\bf d}^m) = \sum_{j=1} W_j(s,{\bf d}^m)\;,
\label{SylvWavesExpand}
\ee
where summation runs over all distinct factors in the set ${\bf d}^m$.
The wave $W_j(s,{\bf d}^m)$ is a quasipolynomial in $s$ 
closely related to prime roots $\rho_j$ of unit.
Namely, Sylvester showed in \cite{Sylv2} that the wave
$W_j(s,{\bf d}^m)$ is a coefficient of
${t}^{-1}$ in the series expansion in ascending powers of $t$ of
\be
F_j(s,t)=\sum_{\rho_j} \frac{\rho_j^{-s} e^{st}}{\prod_{k=1}^{m}
        \left(1-\rho_j^{d_k} e^{-d_k t}\right)}\;.
\label{generatorWj}
\ee
The summation is made over all prime roots of unit
$\rho_j=\exp(2\pi i n/j)$ for $n$ relatively prime to $j$
(including unity) and smaller than $j$.
This result is
just a recipe for calculation of the partition function and it
does not provide explicit formula.

Using the Sylvester recipe we find an explicit formula for the
Sylvester wave $W_j(s,{\bf d}^m)$ in a form of finite sum of the Bernoulli
polynomials of higher
order \cite{bat53,NorlundMemo} multiplied by a periodic function of integer 
period $j$. The periodic factor is expressed through the generalized
Euler polynomials of higher order \cite{Carlitz1960}.
%, which we define similarly to the
%standard Euler polynomials of higher order \cite{bat53,NorlundMemo}.

A special symbolic technique is developed in the theory of polynomials of higher
order, which significantly simplifies computations performed with these 
polynomials. A short description of this technique required
for better understanding of this paper is given in Appendix 
\ref{appendix1}.
%%%%%%%%%%%%%%%%%%%%%%%%%%%%%%%%%%%%%%%%%%%%%%%%%%%%%%%%%%%%%%%%%%%%%
\section{Sylvester wave $W_1(s,{\bf d}^m)$ and Bernoulli polynomials \\
of higher order}
\label{1}
Consider a polynomial part of the partition function
corresponding to the wave $W_1(s,{\bf d}^m)$. It may be found as a residue of
the generator
\be
F_1(s,t) = \frac{e^{st}}{\prod_{i=1}^m (1-e^{-d_i t})}\;.
\label{generatorW1}
\ee
Recalling the generating function for the Bernoulli polynomials of higher 
order \cite{bat53}:
\be
\frac{e^{st} t^m \prod_{i=1}^m d_i}{\prod_{i=1}^m (e^{d_it}-1)} =
\sum_{n=0}^{\infty} B^{(m)}_n(s|{\bf d}^m)
\frac{t^{n}}{n!}\;,
\label{genfuncBernoulli0}
\ee
and a transformation rule
$$
B^{(m)}_n(s|-{\bf d}^m) = B^{(m)}_n(s+\sum_{i=1}^m d_i|{\bf d}^m)\;,
$$
we obtain the relation
\be
\frac{e^{st}}{\prod_{i=1}^m (1-e^{-d_it})} =
\frac{1}{\pi_m} \sum_{n=0}^{\infty} B^{(m)}_n(s+s_m|{\bf d}^m)
\frac{t^{n-m}}{n!}\;,
\label{genfuncBernoulli}
\ee
where
$$
s_m = \sum_{i=1}^m d_i, \ \ \pi_m =  \prod_{i=1}^m d_i\;.
$$

It is immediately seen from (\ref{genfuncBernoulli}) that the coefficient 
of $1/t$ in  (\ref{generatorW1}) is given by the term with $n=m-1$
\begin{equation}
W_1(s,{\bf d}^m) =
\frac{1}{(m-1)!\;\pi_m}
B_{m-1}^{(m)}(s + s_m | {\bf d}^m)\;.
\label{W_1}
\end{equation}
The polynomial part also admits a symbolic form frequently used in theory of
higher order polynomials
\begin{equation}
W_1(s,{\bf d}^m) = \frac{1}{(m-1)!\;\pi_m}
\left(s+s_m + \sum_{i=1}^m d_i \;{}^i\! B\right)^{m-1}\;,
\label{W_1symb}
\end{equation}
where after expansion powers $r_i$ of ${}^i\! B$ are converted into orders
of the Bernoulli numbers
\be
{}^i \! B^{r_i} \Rightarrow B_{r_i}\;.
\label{replacement_rule}
\ee
It is easy to recognize in (\ref{W_1}) the explicit
formula reported recently in
\cite{Beck}, which was obtained by a straightforward computation of the complex 
residue of the generator (\ref{generatorW1}).

Note that basic recursive relation for the Bernoulli polynomials
\cite{NorlundMemo}
\be
B_{n}^{(m)}(s + d_m | {\bf d}^m) -
B_{n}^{(m)}(s | {\bf d}^m) =
n d_m B_{n-1}^{(m-1)}(s | {\bf d}^{m-1})
\label{Bernoulli_recursion}
\ee
naturally leads to the basic recursive relation for the polynomial part of
the partition function:
\be
W_1(s,{\bf d}^m) - W_1(s-d_m,{\bf d}^m) =
W_1(s,{\bf d}^{m-1})\;,
\label{SW1_recursion}
\ee
which coincides with (\ref{SW_recursion}). This indicates that the 
Bernoulli polynomials of higher order represent a natural
basis for expansion of the partition function and its waves.
%%%%%%%%%%%%%%%%%%%%%%%%%%%%%%%%%%%%%%%%%%%%%%%%%%%%%%%%%%%%%%%%%%%%%
\section{Sylvester wave $W_2(s,{\bf d}^m)$ and Euler numbers \\
of higher order}
\label{2}
In order to compute the Sylvester wave with period $j>1$ we note, that
the summand in the expression (\ref{generatorWj}) can be rewritten as
a product
\be
F_j(s,t) = \sum_{\rho_j}
\frac{e^{st}}{\prod_{i=1}^{\omega_j} (1-e^{-d_it})} \times
\frac{\rho_j^{-s}}{\prod_{i=\omega_j+1}^m (1-\rho_j^{d_i}e^{-d_i t})}\;,
\label{generator_product}
\ee
where the elements in ${\bf d}^m$ are sorted in a way that $j$ is a divisor
for first $\omega_j$ elements (we say that $j$ has weight
$\omega_j$), and the rest elements in the set are not
divisible by $j$.

Then a 2-periodic
Sylvester wave $W_2(s,{\bf d}^m)$ is a residue of the generator
\be
F_2(s,t) =
\frac{e^{st}}{\prod_{i=1}^{\omega_2} (1-e^{-d_it})} \times
\frac{(-1)^{s}}{\prod_{i=\omega_2+1}^m (1+e^{-d_i t})}\;,
\label{generatorW2}
\ee
where first $\omega_2$ integers $d_i$ are even, and the summation is 
omitted being trivially restricted to the only value $\rho_2=-1$.
Recalling the generating function for the Euler polynomials of higher order
\cite{bat53, NorlundMemo} and corresponding recursive relation 
\begin{eqnarray}
\frac{2^m e^{st}}{\prod_{i=1}^m (e^{d_i t}+1)} =
\sum_{n=0}^{\infty} E_n^{(m)}(s | {\bf d}^m) \frac{t^n}{n!}\;,
\label{Euler_GF} \\
E_n^{(m)}(s+d_m | {\bf d}^m)+E_n^{(m)}(s | {\bf d}^m)=
2E_n^{(m-1)}(s | {\bf d}^{m-1})\;,\nonumber
\end{eqnarray}
we may rewrite (\ref{generatorW2}) as double infinite sum
\be
\frac{(-1)^{s}}{2^{m-\omega_2} \pi_{\omega_2}}
\sum_{n=0}^{\infty} B^{(\omega_2)}_n(s+s_{\omega_2}|{\bf d}^{\omega_2})
\frac{t^{n-\omega_2}}{n!}
\sum_{l=0}^{\infty} E_l^{(m-\omega_2)}(s_m-s_{\omega_2} | {\bf d}^{m-\omega_2})
\frac{t^l}{l!}\;.
\label{g2sum}
\ee
The coefficient of $1/t$ in the above series is found for
$n+l=\omega_2-1$, so that we end up with a finite sum:
\be
W_2(s,{\bf d}^m) =
\frac{(-1)^{s}}{(\omega_2-1)!\; 2^{m-\omega_2} \pi_{\omega_2}}
\sum_{n=0}^{\omega_2-1} \binom{\omega_2-1}{n}
B^{(\omega_2)}_n(s+s_{\omega_2}|{\bf d}^{\omega_2})
E_{\omega_2-1-n}^{(m-\omega_2)}(s_m-s_{\omega_2} | {\bf d}^{m-\omega_2}). %,
\label{W2}
\ee
This expression may be rewritten as a symbolic power similar to
(\ref{W_1symb}):
\be
W_2(s,{\bf d}^m) =
\frac{(-1)^{s}}{(\omega_2-1)! \; 2^{m-\omega_2} \pi_{\omega_2}}
\left(
s+s_m + \sum_{i=1}^{\omega_2} d_i \;{}^i\! B +
\sum_{i=\omega_2+1}^{m} d_i \;{}^i\! E(0)
\right)^{\omega_2-1},
\label{W2symb}
\ee
where the rule for the Euler polynomials at zero $E_n(0)$ similar to
(\ref{replacement_rule}) is applied.
It is easy to rewrite formula (\ref{W2symb}) in a form
\be
W_2(s,{\bf d}^m) =
\frac{(-1)^{s}}{(\omega_2-1)!\; 2^{m-\omega_2} \pi_{\omega_2}}
\sum_{n=0}^{\omega_2-1} \binom{\omega_2-1}{n}
B^{(\omega_2)}_n(s+s_{m}|{\bf d}^{\omega_2})
E_{\omega_2-1-n}^{(m-\omega_2)}(0|{\bf d}^{m-\omega_2}),
\label{W2last}
\ee
where $E_{n}^{(m)}(0|{\bf d}^{m})$ denote the Euler polynomials
of higher orders computed at zero as follows:
\be
E_{n}^{(m)}(0|{\bf d}^{m}) =
\left[
\sum_{i=1}^{m} d_i \;{}^i\! E(0)
\right]^{n}.
\label{Enumbers}
\ee
The formula (\ref{W2last}) shows that the wave $W_2(s,{\bf d}^{m})$ can be
written as an expansion
over the Bernoulli polynomials of higher order with constant coefficients,
multiplied by a 2-periodic function $(-1)^s$.
%%%%%%%%%%%%%%%%%%%%%%%%%%%%%%%%%%%%%%%%%%%%%%%%%%%%%%%%%%%%%%%%%%%%%
\section{Sylvester waves $W_j(s,{\bf d}^m) \ (j>2)$ and Euler \\
polynomials of higher order}
\label{j}
Frobenius \cite{Frobenius} studied in great detail the polynomials
$H_n(s,\rho)$ satisfying the generating function 
%(see \cite{Carlitz59} for further reference)
\be
\frac{(1-\rho) e^{st}}{e^t-\rho} = \sum_{n=0}^{\infty} H_n(s,\rho)
\frac{t^n}{n!}, \ \ (\rho \ne 1),
\label{EulerRegGFnew}
\ee
which reduces to definition of the Euler polynomials at fixed value of the
parameter $\rho$
$$
E_n(s) = H_n(s,-1).
$$
The polynomials $H_n(\rho) \equiv H_n(0,\rho)$ satisfy the
symbolic recursion ($H_0(\rho)=1$)
\be
\rho H_n(\rho) = (H(\rho)+1)^n, \ \ \ n>0.
\label{EulerSymbolic}
\ee
The generalization of (\ref{Euler_GF}) is straightforward
\be
\frac{ e^{st} \prod_{i=1}^m (1-\rho^{d_i})}
{\prod_{i=1}^m (e^{d_i t}-\rho^{d_i})} =
\sum_{n=0}^{\infty} H_n^{(m)}(s, \rho | {\bf d}^m) \frac{t^n}{n!},
\ \ (\rho^{d_i} \ne 1),
\label{Euler_GFnew}
\ee
where the corresponding recursive relation for $H_n^{(m)}(s, \rho | {\bf 
d}^m)$ has the form
\begin{equation}
H_n^{(m)}(s+d_m,\rho | {\bf d}^m)-\rho^{d_m}H_n^{(m)}(s,\rho | {\bf d}^m)=
\left(1-\rho^{d_m}\right)H_n^{(m-1)}(s,\rho | {\bf d}^{m-1})\;.
\label{Euler_GFnewrecur}
\end{equation}
The {\em generalized Euler polynomials of higher order}
$H_n^{(m)}(s, \rho | {\bf d}^m)$
introduced by L. Carlitz in \cite{Carlitz1960}
can be defined through the symbolic
formula
\be
H_n^{(m)}(s, \rho | {\bf d}^m) =
\left(s + \sum_{i=1}^{m} d_i \;{}^i\! H(\rho^{d_i}) \right)^n,
\label{EulerNewSymbolic}
\ee
where $H_n(\rho)$ computed from the
relation
$$
\frac{1-\rho}{e^t-\rho} = \sum_{n=0}^{\infty} H_n(\rho)
\frac{t^n}{n!},
$$
or using the recursion (\ref{EulerSymbolic}). Using the polynomials 
$H_n^{(m)}(s, \rho | {\bf d}^m)$ we can compute Sylvester wave of 
arbitrary period.

Consider a $j$-periodic Sylvester wave $W_j(s,{\bf d}^m)$,
and rewrite the summand in (\ref{generator_product}) as double infinite sum
\be
\frac{\rho_j^{-s}}{\pi_{\omega_j} \; \prod_{i=\omega_j+1}^m (1-\rho_j^{d_i})}
\sum_{n=0}^{\infty} B^{(\omega_j)}_n(s+s_{\omega_j}|{\bf d}^{\omega_j})
\frac{t^{n-\omega_j}}{n!}
\sum_{l=0}^{\infty} H_l^{(m-\omega_j)}(s_m-s_{\omega_j}, \rho_j | {\bf
 d}^{m-\omega_j})
\frac{t^l}{l!}.
\label{gjsum}
\ee
The coefficient of $1/t$ in the above series is found for
$n+l=\omega_j-1$, so that we have a finite sum:
\bea
W_j(s,{\bf d}^m) & = &
\frac{1}{(\omega_j-1)! \; \pi_{\omega_j}}
\sum_{\rho_j}
\frac{\rho_j^{-s}}{\prod_{i=\omega_j+1}^m (1-\rho_j^{d_i})}
\times \nonumber \\
&&\sum_{n=0}^{\omega_j-1} \binom{\omega_j-1}{n}
B^{(\omega_j)}_n(s+s_{\omega_j}|{\bf d}^{\omega_j})
H_{\omega_j-1-n}^{(m-\omega_j)}(s_m-s_{\omega_j}, \rho_j | {\bf
 d}^{m-\omega_j})\;. 
\label{Wj}
\eea
This expression may be rewritten as a symbolic power similar to
(\ref{W2symb}):
\be
W_j(s,{\bf d}^m) =
\frac{1}{(\omega_j-1)! \; \pi_{\omega_j}}
\sum_{\rho_j}
\frac{\rho_j^{-s}}{\prod_{i=\omega_j+1}^m (1-\rho_j^{d_i})}
\left(
s+s_m + \sum_{i=1}^{\omega_j} d_i \;{}^i\! B + \!\!\!
\sum_{i=\omega_j+1}^{m} \!\!\! d_i \;{}^i\! H(\rho_j^{d_i})
\right)^{\omega_j-1} \!\!\!\!\!\!\!\!\;,
\label{Wjsymb}
\ee
which is equal to
\bea
W_j(s,{\bf d}^m) & = &
\frac{1}{(\omega_j-1)! \; \pi_{\omega_j}}
\sum_{n=0}^{\omega_j-1} \binom{\omega_j-1}{n}
B^{(\omega_j)}_n(s+s_{m}|{\bf d}^{\omega_j})
\times \nonumber \\
&& \sum_{\rho_j}
\frac{\rho_j^{-s}}{\prod_{i=\omega_j+1}^m (1-\rho_j^{d_i})}
H_{\omega_j-1-n}^{(m-\omega_j)}[\rho_j |{\bf d}^{m-\omega_j}]\;,
\label{WjBern}
\eea
where
\be
H_{n}^{(m)}[\rho |{\bf d}^{m}] =
H_{n}^{(m)}(0,\rho |{\bf d}^{m}) =
\left[\sum_{i=1}^{m} d_i \;{}^i\! H(\rho^{d_i})\right]^n\;,
\label{Zrhonumbers}
\ee
are generalized Euler numbers of higher order
and it is assumed that
$$
H_{0}^{(0)}[\rho | \emptyset] = 1, \
H_{n}^{(0)}[\rho | \emptyset] = 0\;, \ n>0\;.
$$

It should be underlined that the presentation of the Sylvester wave as a
finite sum of the Bernoulli polynomials with periodic coefficients
(\ref{WjBern}) is not unique. The symbolic formula (\ref{Wjsymb}) can be 
cast into a sum of the generalized Euler polynomials
\bea
W_j(s,{\bf d}^m) & = &
\frac{1}{(\omega_j-1)! \; \pi_{\omega_j}}
\sum_{n=0}^{\omega_j-1} \binom{\omega_j-1}{n}
B^{(\omega_j)}_n[{\bf d}^{\omega_j}]
\times \nonumber \\
&& \sum_{\rho_j}
\frac{\rho_j^{-s}}{\prod_{i=\omega_j+1}^m (1-\rho_j^{d_i})}
H_{\omega_j-1-n}^{(m-\omega_j)}(s+s_m, \rho_j |{\bf d}^{m-\omega_j})\;,
\label{WjEuler}
\eea
where 
$$
B^{(m)}_n[{\bf d}^{m}] = B^{(m)}_n(0|{\bf d}^{m})
$$
are the Bernoulli numbers of higher order.

Substitution of the expression (\ref{WjBern}) into the expansion 
(\ref{SylvWavesExpand}) immediately produces the partition function 
$W(s,{\bf d}^{m})$ as finite sum of the Bernoulli
polynomials of higher order multiplied by periodic functions with
period equal to the least common multiple of the elements in ${\bf d}^m$
\bea
W(s,{\bf d}^m) & = & \sum_j
\frac{1}{(\omega_j-1)! \; \pi_{\omega_j}}
\sum_{n=0}^{\omega_j-1} \binom{\omega_j-1}{n}
B^{(\omega_j)}_n(s+s_{m}|{\bf d}^{\omega_j})
\times \nonumber \\
&& \sum_{\rho_j}
\frac{\rho_j^{-s}}{\prod_{i=\omega_j+1}^m (1-\rho_j^{d_i})}
H_{\omega_j-1-n}^{(m-\omega_j)}[\rho_j | {\bf d}^{m-\omega_j}]\;.
\label{WBern}
\eea

The partition function $W(s,{\bf d}^m)$ has several interesting properties.
Analysis of the generating function (\ref{genfunc}) shows that the partition
function is a homogeneous function of zero order with respect to all 
its arguments, i.e.,
\be
W(k s,k {\bf d}^m) = W(s,{\bf d}^m).
\label{homogeneity}
\ee
This property appears very useful for computation of the partition
function in case when the elements $d_i$ have a common factor $k$, then
\be
W(s, k {\bf d}^m) = W\left(\frac{s}{k},{\bf d}^m\right).
\label{factor}
\ee

The case of $m$ identical elements ${\bf p}^m=\{p,\ldots,p\}$ appears to be the
simplest and is reduced to the known formula for Catalan partitions 
\cite{catal838}: {\it the Diophantine equation  $x_1+x_2+\;.\;.\;.\;+x_m=s$ has 
${s+m-1\choose s}$ sets of non-negative solutions.}
%, the repetitions are not excluded.}

Using (\ref{factor}) for $s$ divisible by $p$ we arrive at
$$
W(s,{\bf p}^m) = W\left(\frac{s}{p},{\bf 1}^m\right) =
W_1\left(\frac{s}{p},{\bf 1}^m\right) =
\frac{B_{m-1}^{(m)}(s/p + m | {\bf 1}^m)}{(m-1)!}.
$$
The straightforward computation shows that
$$
B_{m-1}^{(m)}(s + m | {\bf 1}^m) = \prod_{k=1}^{m-1} (s+k) =
\frac{(s+m-1)!}{s!},
$$
so that
\be
W(s,{\bf p}^m) =
\left\{ \begin{array}{ll}
\prod_{k=1}^{m-1} \left(1+\frac{s}{kp} \right), & s=0 \pmod p,\\
0 \;, & s \ne 0  \pmod p.
\end{array}\right.
\label{identical_d}
\ee
In the end of this Section we consider a special case of the tuple
$\{p_1,p_2,\ldots p_m\}$ of primes $p_j$ which leads to essential 
simplification of the formula (\ref{WBern}). 
The first Sylvester wave $W_1$ is given by (\ref{W_1}) while
all higher waves arising
% from (\ref{WjEuler}) for all $\omega_r=1,\;r\geq 2$
are purely periodic
\begin{equation}
W_{p_{i}}(s;\{p_1,p_2,\ldots,p_m\})=
%\frac{1}{p_i}
%\sum_{k=1\atop k\neq 0 \bmod p_{i}}^{k=p_{i}-1}\rho_{p_{i}}^{-ks}
%\sum_{k=1}^{k=\;p_{i}-1}\rho_{p_{i}}^{-ks}
%\prod_{j=2\atop p_{j}\neq 0 \bmod p_{i}}^m
%\prod_{j=1 \atop j \ne i}^m
%\frac{1}{1-\rho_{p_{i}}^{kp_{j}}}=
\frac{1}{p_i}\sum_{k=1}^{p_{i}-1}\frac{\rho_{p_{i}}^{-ks}}
{\prod_{j\neq i}^m\left(1-\rho_{p_{i}}^{kp_{j}}\right)}\;.
\label{prim2}
\end{equation}
The further simplification $m=2,\,s=ap_1p_2$ makes it possible to
verify the partition identity
\begin{eqnarray}
W(a p_1p_2,\{p_1,p_2\})=a+1\;,
\label{aa1}
\end{eqnarray}
which follows from the recursion relation (\ref{SW_recursion}) for 
the restricted partition function and its definition 
\begin{eqnarray}
&&W(ap_1p_2,\{p_1,p_2\}) - W(ap_1p_2-p_1,\{p_1,p_2\}) =
W(ap_1p_2,\{p_2\})\;,\nonumber\\
&&W(ap_1p_2,\{p_2\})=1\;,\;\;\;W(ap_1p_2-p_1,\{p_1,p_2\})=
W((a-l)p_1p_2+(lp_2-1)p_1,\{p_1,p_2\})=a\;,\nonumber
\end{eqnarray}
where $a$ solutions of the Diophantine equation
$p_1X+p_2Y=(a-l)p_1p_2+(lp_2-1)p_1$ correspond to $l=1,\ldots,a$. 
The relation (\ref{aa1}) has an important geometrical interpretation, namely,
a line $p_1X+p_2Y=ap_1p_2$ in the $XY$ plane
passes exactly through $a+1$ points with non-negative integer
coordinates.

The verification of (\ref{aa1}) is straightforward
 %Indeed, performing summation in (\ref{prim2}) we obtain 
(see Appendix B for details):
\begin{eqnarray}
&&W_1(ap_1p_2,\{p_1,p_2\})=a+\frac{1}{2}\left(\frac{1}{p_1}+
\frac{1}{p_2}\right)\;,\nonumber\\
&&W_{p_1}(ap_1p_2,\{p_1,p_2\})=\frac{1}{2}-\frac{1}{2p_1}\;,\;\;\;
W_{p_2}(ap_1p_2,\{p_1,p_2\})=\frac{1}{2}-\frac{1}{2p_2}\;,
\label{sylvp1p2}
\end{eqnarray}
which produces the required result.

A generalization of (\ref{aa1}) is possible using the explicit form of
the partition function
\be
W(s,\{p_1,p_2\}) =
\frac{1}{p_1p_2}\left(s+\frac{p_1+p_2}{2}\right)+
\frac{1}{p_1} \sum_{\rho_{p_1}} \frac{\rho_{p_1}^{-s}}{1-\rho_{p_1}^{p_2}} +
\frac{1}{p_2} \sum_{\rho_{p_2}} \frac{\rho_{p_2}^{-s}}{1-\rho_{p_2}^{p_1}}.
\label{2primes}
\ee
Setting here $s=ap_1p_2+b, \, 0 \le b < p_1p_2$ and
noting that the value of two last terms in (\ref{2primes})
don't depend on the integer $a$, one can easily see that
\be
W(ap_1p_2+b,\{p_1,p_2\}) = a + W(b,\{p_1,p_2\}),
\label{reduct}
\ee
which reduces the procedure to computation of the first $p_1p_2$ 
values of $W(s,\{p_1,p_2\})$. Recalling that $W(0,\{p_1,p_2\}) = 1$ we 
immediately recover (\ref{aa1}) as a particular case of (\ref{reduct}).

%%%%%%%%%%%%%%%%%%%%%%%%%%%%%%%%%%%%%%%%%%%%%%%%%%%%%%%%%%%%%%%%%%%%%
\section{Recursive Relation for Sylvester Waves}
\label{recurs}
In this Section we prove  that the recursive relation similar to 
(\ref{SW_recursion}) holds not only for the entire partition function 
$W(s,{\bf d}^m)$ and its polynomial part $W_1(s,{\bf d}^m)$ but also
for each Sylvester wave
\be
W_j(s,{\bf d}^m) - W_j(s-d_m,{\bf d}^m) =
W_j(s,{\bf d}^{m-1})\;.
\label{SWj_recursion}
\ee

When $j$ is not a divisor of $d_m$, the weight $\omega_j$ doesn't change
in transition from ${\bf d}^{m-1}$ to ${\bf d}^{m}$. Denoting for brevity
$$
A(s) = s+s_{m-1} + \sum_{i=1}^{\omega_j}d_i \;{}^i\! B + \!\!\!
\sum_{i=\omega_j+1}^{m-1} \!\! d_i \;{}^i\! H(\rho_j^{d_i})\;, \ \
B_{\omega_j} = \frac{1}{(\omega_j-1)! \; \pi_{\omega_j}}\;,
$$
we have
\bea
W_j(s,{\bf d}^m)  & = &
B_{\omega_j} %\frac{1}{(\omega_j-1)! \; \pi_{\omega_j}}
\sum_{\rho_j}
\frac{\rho_j^{-s}}{\prod_{i=\omega_j+1}^m (1-\rho_j^{d_i})}
\left(
A(s) + d_m[1 + H(\rho_j^{d_m})]
\right)^{\omega_j-1} \nonumber \\
&=&
B_{\omega_j}
\sum_{\rho_j}
\frac{\rho_j^{-s}}{\prod_{i=\omega_j+1}^m (1-\rho_j^{d_i})}
\sum_{l=0}^{\omega_j-1}
\binom{\omega_j-1}{l} A^{\omega_j-1-l}(s) d_m^l
[1 + H(\rho_j^{d_m})]^l\;. \nonumber
\eea
Now using (\ref{EulerSymbolic}) we have
\bea
W_j(s,{\bf d}^m)  & = &
B_{\omega_j}
\sum_{\rho_j}
\frac{\rho_j^{-s}}{\prod_{i=\omega_j+1}^m (1-\rho_j^{d_i})}
\left\{A^{\omega_j-1}(s) + \rho_j^{d_m} \sum_{l=1}^{\omega_j-1}
\binom{\omega_j-1}{l} A^{\omega_j-1-l}(s) d_m^l H_l(\rho_j^{d_m})
\right\} \nonumber \\
&=&
B_{\omega_j}
\sum_{\rho_j}\frac{\rho_j^{-s}}{\prod_{i=\omega_j+1}^m (1-\rho_j^{d_i})}
\left\{(1-\rho_j^{d_m})A^{\omega_j-1}(s) + \rho_j^{d_m}\left(
A(s) + d_m H(\rho_j^{d_m})\right)^{\omega_j-1}\right\} \nonumber \\
&=&
B_{\omega_j}
\sum_{\rho_j}
\frac{\rho_j^{-(s-d_m)}}{\prod_{i=\omega_j+1}^m (1-\rho_j^{d_i})}
\left(
A(s) + d_m H(\rho_j^{d_m})
\right)^{\omega_j-1}
+ 
%\\ &&
B_{\omega_j}
\sum_{\rho_j}
\frac{\rho_j^{-s}A^{\omega_j-1}(s)}
{\prod_{i=\omega_j+1}^{m-1} (1-\rho_j^{d_i})}
\nonumber \\
&=&
W_j(s-d_m,{\bf d}^m) +
W_j(s,{\bf d}^{m-1})\;. 
\eea
In case of $j$ being divisor of $d_m$ the weight of $j$ 
for the set ${\bf d}^{m-1}$
is equal to $\omega_j-1$, and we have
\be
W_j(s,{\bf d}^{m-1}) =
\frac{(\omega_j-1) d_m}{(\omega_j-1)! \; \pi_{\omega_j}}
\sum_{\rho_j}
\frac{\rho_j^{-s}}{\prod_{i=\omega_j}^{m-1} (1-\rho_j^{d_i})}
\left(
s+s_{m-1} + \sum_{i=1}^{\omega_j-1} d_i \;{}^i\! B + \!\!\!
\sum_{i=\omega_j}^{m-1} \!\!\! d_i \;{}^i\! H(\rho_j^{d_i})
\right)^{\omega_j-2} \!\!\!\!\!\!\!\!\;.
\label{Wm-1}
\ee
Denoting
$$
A(s) = s+s_{m-1} + \sum_{i=1}^{\omega_j-1} d_i \;{}^i\! B + \!\!
\sum_{i=\omega_j}^{m-1} \! d_i \;{}^i\! H(\rho_j^{d_i}), \ \ \
D(s,\rho_j) = \frac{\rho_j^{-s}}{\prod_{i=\omega_j+1}^m (1-\rho_j^{d_i})}\;,
$$
and using the symbolic formula for the Bernoulli numbers
\cite{NorlundMemo}
$$
(B+1)^n = B^n = B_n \ \ (n \ne 1)\;,
$$
we obtain
\bea
W_j(s,{\bf d}^m)  & = &
B_{\omega_j} \sum_{\rho_j}
D(s,\rho_j)
[A(s) + d_m(B+1)]^{\omega_j-1} \nonumber \\
&=&
B_{\omega_j}
\sum_{\rho_j}
D(s,\rho_j)
\sum_{l=0}^{\omega_j-1}
\binom{\omega_j-1}{l} A^{\omega_j-1-l}(s) d_m^l
(B+1)^l \\
&=&
B_{\omega_j}
\sum_{\rho_j}
D(s,\rho_j)
[A(s) + d_m B]^{\omega_j-1}
+
B_{\omega_j} d_m (\omega_j-1)
\sum_{\rho_j}
D(s,\rho_j) A^{\omega_j-2}(s)
 \nonumber \\
&=&
W_j(s-d_m,{\bf d}^m) +W_j(s,{\bf d}^{m-1})\;, \nonumber
\eea
which completes the proof.
%%%%%%%%%%%%%%%%%%%%%%%%%%%%%%%%%%%%%%%%%%%%%%%%%%%%%%%%%%%%%%%%%%%%%
\section{Partition function $W\left(s,\{\overline{m}\}\right)$ for a set 
of natural \\numbers \label{Sm}}

Sylvester waves for a set of consecutive natural numbers 
$\{1,2,\dots,m\}=\{\overline{m}\}$ was under special consideration in 
\cite{Rama}. An importance of this case based on its relation to the
invariants of symmetric group $S_m$ (see next Section) and, second, 
$W(s,\{\overline{m}\})$ form a natural basis to utilize the 
partition functions for every subsets of $\{1,2,\dots,m\}$. This case is 
also important due to the famous Rademacher formula \cite{Radem37} for 
{\it unrestricted partition function} $W(s,\{\overline{s}\})$, but 
the latter already belongs to the analytical number theory.

The representation for $W(s,\{\overline{m}\})$ in terms of higher 
Bernoulli polynomials comes when we put into (\ref{WBern})
\begin{equation}
\omega_j=\left[\frac{m}{j}\right]\;,\;\;\;\;
\pi_{\omega_j}=\omega_j!\;j^{\omega_j}\;,\;\;\;\;
s_{\omega_j}=\frac{\omega_j(\omega_j+1)}{2}\;,
\label{symmet}
\end{equation}
where $[x]$ denotes integer part of $x$. The partition function in this 
case reads
\bea
W(s,\{\overline{m}\}) & = & \sum_{j=1}^m
\frac{j^{-\omega_j}}{(\omega_j-1)!\; \omega_j!}
\sum_{n=0}^{\omega_j-1} \binom{\omega_j-1}{n}
B^{(\omega_j)}_n\left(s+\frac{m(m+1)}{2}|{\bf d}^{\omega_j}\right)
\times \nonumber \\
&& \sum_{\rho_j} \frac{\rho_j^{-s}}
{\prod_{i=\omega_j+1}^m (1-\rho_j^{d_i})}
H_{\omega_j-1-n}^{(m-\omega_j)}[\rho_j|{\bf d}^{m-\omega_j}]\;,
\label{WSmWODenom}
\eea
where ${\bf d}^{\omega_j} = j\{\overline{\omega}\}$, so that for each $j$
we have elements
divisible by $j$ at first $\omega_j$ positions. The expression for the 
Sylvester wave of the maximal period $m$ looks particularly simple
\be
W_m(s,\{\overline{m}\}) = \frac{1}{m^2} \sum_{\rho_m}
\rho_m^{-s}\;.
\label{WmSm}
\ee
The straightforward calculations show that the expression (\ref{WSmWODenom}) 
produces exactly the same formulas for $m=1,2,\ldots,12$ which were obtained 
in \cite{Rama}.
\begin{figure}[t]
\psfig{figure=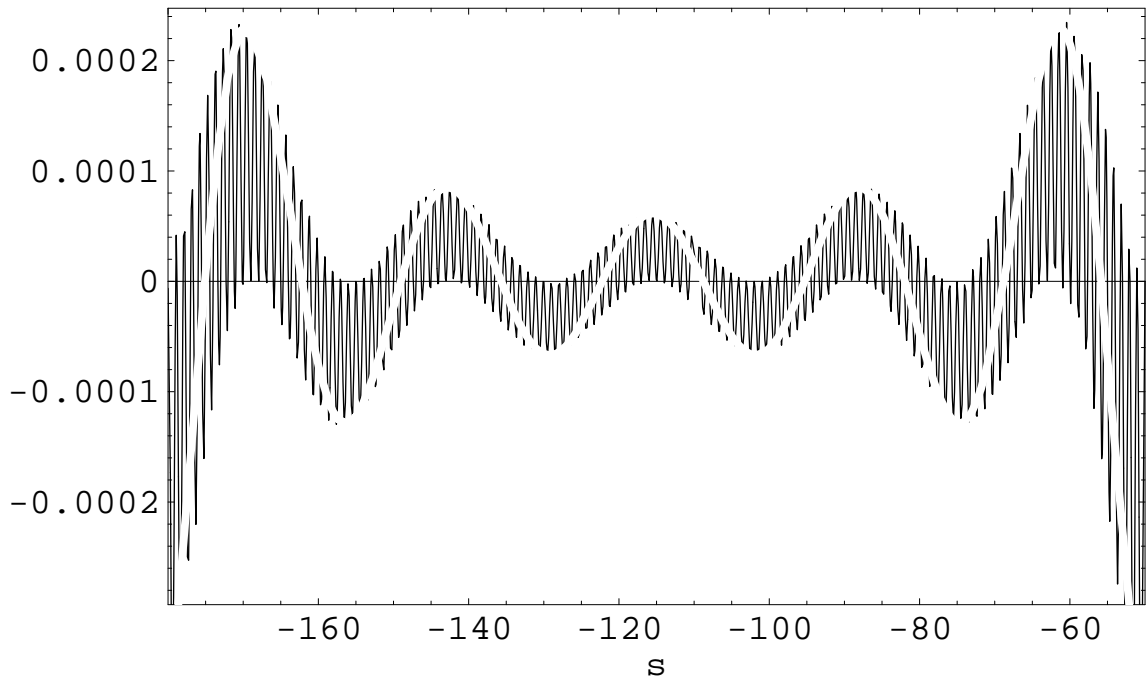,width=6in}
\caption{Plots of the partition function $W(s,\{\overline{21}\})$ ({\it
black curve}) and its first Sylvester wave $W_1(s,\{\overline{21}\})$
({\it white curve}) showing that the polynomial part provides an
important information about the partition function behavior.}
\label{W21approx}
\end{figure}

It needs to be noted that typically the argument $s$ in all formulas
derived above is assumed to have integer values, but it is obvious that all
results can be extended to real values of $s$, though such extension is not
unique.
Continuous values of the argument provide a convenient way to analyze the
behavior of the partition function and its waves. In this work we choose
the natural extension scheme based on the trigonometric functions
$$
\rho_j^s = e^{2 \pi i n s/j} = \cos \frac{2 \pi n s}{j} +
i \sin \frac{2 \pi n s}{j}.
$$
We finish this Section with a brief discussion of a phenomenon better
observed in graphics of $W(s,\{\overline{m}\})$
with large $m$ rather from the explicit expressions (see formulas (52) and
Figures of restricted partition functions in \cite{Rama}).

\begin{figure}[t]
\psfig{figure=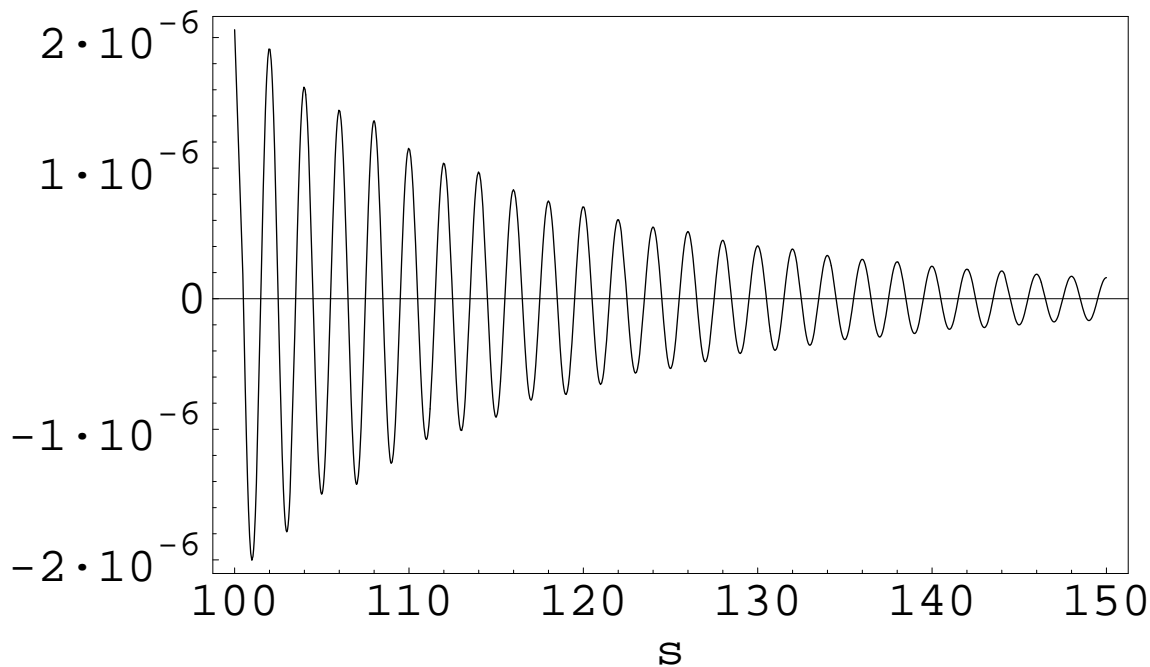,width=6in}
\caption{Plot of the normalized difference
$[W(s,\{\overline{21}\})/W_1(s,\{\overline{21}\})-1]$
showing that the polynomial
part $W_1(s,\{\overline{21}\})$ at large values of the argument $s$ gives
a very accurate approximation to the partition function
$W(s,\{\overline{21}\})$.}
\label{W21diff}
\end{figure}

In the range $[-\frac{m(m+1)}{2},0]$ where $W(s,\{\overline{m}\})$ has
all its zeroes, one can easily assume an existence
of a function $\widetilde{W}(s,\{\overline{m}\})$ which envelopes
$W(s,\{\overline{m}\})$ or approximates it in some sense. The
decomposition of $W(s,\{\overline{m}\})$ into the Sylvester waves shows
that this role may be assigned to the wave $W_1(s,\{\overline{m}\})$. The
Figures \ref{W21approx}, \ref{W21diff} show that $W_1(s,\{\overline{21}\})$
serves as a good approximant for $W(s,\{\overline{21}\})$ in this range
as well as for large $s$.
%%%%%%%%%%%%%%%%%%%%%%%%%%%%%%%%%%%%%%%%%%%%%%%%%%%%%%%%%%%%%%%%%%%%%
\section{Application to invariants of finite groups}
\label{finitegroup}
The restricted partition function $W(s,{\bf d}^m)$ has a strong
relationship to the invariants of finite reflection groups $G$
acting on the vector space $V$ over the field of complex numbers. If $M^G(t)$
is a Molien function of the finite group, $d_{r}$ and $m$ are degrees
and a number of the basic homogeneous invariants respectively, then its
series expansion
% (\ref{molien1}) 
in $t$ gives a number $P(s,G)$ of
algebraically independent invariants of the degree $s$. The set of natural 
numbers $\{\overline{m}\}$ corresponds to the symmetric group
$S_m:\;W(s,\{\overline{m}\})=P(s,S_m)$. The list of $P(s,G)$ for all 
indecomposable reflections groups $G$ acting over the field of real 
numbers and known as {\it Coxeter groups} is presented in \cite{Rama}. It 
is easy to extend these formulas over indecomposable pseudoreflections 
groups acting over the field of complex numbers using the list of 37 
groups given by Shepard and Todd \cite{Shepar54}. In this Section we 
extend the results of Section \ref{j} to all finite groups.

First, we recall an algebraic setup of the problem.
The fundamental problem of the invariant theory consists in determination of an
algebra ${\sf R}^G$ of invariants. Its solution is given by the
Noether theorem \cite{Benson93}: ${\sf R}^G$ is generated by a
polynomial $\vartheta_k(x_j)$ as an algebra due to action of finite group
$G\subset GL(V^q)$ on the $q$-dimensional vector space $V^q(x_j)$
over the field of complex numbers by not more than ${|G|+q\choose q}$
homogeneous invariants, of degrees not exceeding the order $|G|$ of group
\be
k\leq {|G|+q\choose q}\;,\;\;\;j\leq \dim V^q=q\;,\;\;\;
\deg \vartheta_k(x_j)\leq |G|\;.
\label{molien0}
\ee
To enumerate the invariants explicitly, it is convenient to classify
them by their degrees (as polynomials). A classical theorem of
Molien \cite{Benson93} gives an explicit expression for a number $P(s,G)$
of all homogeneous invariants of degree $s$
\be
M^G(t)=\frac{1}{|G|}\sum_{l=1}^{|G|}
\frac{\widetilde{\chi}({\widehat g}_l)}
{\det({\hat I} -t\; {\widehat g}_l)}=
\sum_{s=0}^{\infty} P(s,G) t^s\;,\;\;\;\;P(0,G)=1\;,
\label{molien1}
\ee
where ${\widehat g}_l$ are non--singular $(n\times n)$--permutation
matrices with entries, which form the regular representation of $G$,
${\widehat I}$ is the identity matrix and $\widetilde{\chi}$ is the
complex conjugate to character $\chi$. The further progress is due to
Hilbert and his {\it syzygy theorem} \cite{Benson93}. For our purpose
it is important that $M^G(t)$ is a rational polynomial
\be
M^G(t)=\frac{N^G(t)}{\prod_{l=1}^n \left(1-t^{d_l}\right)}\;,\;\;\;\;
N^G(t)=\sum_{k=0}Q(k,G)\;t^k.
\label{molien2a}
\ee
The formula (\ref{molien2a}) is very convenient to express the function
$P(s,G)$ through the Sylvester waves $W(s,{\bf d}^m)$. Recalling the
definition (\ref{genfunc}) of the generating function $F(t,{\bf d}^m)$
consider a general term $t^k F(t,{\bf d}^m)$ of the Molien function
(\ref{molien2a})
\bea
t^k F(t,{\bf d}^m) = \sum_{s=0}^{\infty} W(s,{\bf d}^m) t^{s+k} =
\sum_{s=k}^{\infty} W(s-k,{\bf d}^m) t^{s},
\label{gen_term_Molien}
\eea
so that the corresponding partition function is $W(s-k,{\bf d}^m)$, which 
implies that the number $P(s,G)$ of all homogeneous invariants of degree 
$s$ for the finite group $G$ can be expressed through the simple relation
\be
P(s,G)=\sum_{k=0}^{s}Q(k,G) W(s-k,{\bf d}^m)\;.
\label{fin2}
\ee

We consider several instructive examples for which the explicit
expression of the Molien function $M^G(t)$ and the corresponding number of
homogeneous invariants $P(s,G)$ are given.

1. Alternating group ${\sf A}_n$ generated by its natural $n$--dimensional
representation, $|{\sf A}_n|=n!/2$.
\bea
M_{{\sf A}_n}(t) & = & \left[1+t^{\binom{n}{2}}\right]
\prod_{k=1}^n\frac{1}{1-t^k}\;. \nonumber \\
P(s,{\sf A}_n) & = &
W(s,\{\overline{n}\}) + W\left(s-\frac{n(n-1)}{2},\{\overline{n}\}\right).
\label{altern}
\eea
The group ${\sf A}_n$ is acting on Euclidean vector space ${\mathbb R}^n$.

2. Group ${\sf G}_2$ generated by matrix {\footnotesize $
\left(\begin{array}{cc}
\rho_{n} & 0 \\
0 & \rho_{n}^{-1}\end{array}\right)$}, where $\rho_{n}=e^{2\pi i/n}$ is
a primitive $n$--th root of unity, $|{\sf G}_2|=n$.
\bea
M_{{\sf G}_2}(t)&=&
%\frac{1}{n}\sum_{k=1}^n\frac{1}{(1-\rho_{n,k} t)(1-\rho_{n,k}^{-1} t)}=
\frac{1+t^n}{(1-t^2)(1-t^n)}\;. \nonumber \\
P(s,{\sf G}_2) & = &
W(s,\{2,n\}) +  W(s-n,\{2,n\}).
\label{rotation}
\eea
${\sf G}_2$ is isomorphic as an abstract group to the cyclic group ${\sf 
Z}_n$ acting on Euclidean vector space ${\mathbb R}^2$.

3. Group ${\sf G}_3$ generated by the matrices {\footnotesize
$\left(\begin{array}{cc}
\rho_{n} & 0 \\
0 & \rho_{n}^{-1}\end{array}\right)$} and {\footnotesize
$\left(\begin{array}{cc}
0 & 1 \\
1 & 0\end{array}\right)$}, $|{\sf G}_3|=2n$.
\bea
M_{{\sf G}_3}(t)=\frac{1}{(1-t^2)(1-t^n)}\;,\;\;\;
P(s,{\sf G}_3)=W(s,\{2,n\})
\label{dihedr}
\eea
${\sf G}_3$ is isomorphic as an abstract group to the dihedral group ${\sf
I}_n$ acting on Euclidean vector space ${\mathbb R}^2$.

4. Group ${\sf G}_4$ generated by $(n\times n)$--diagonal matrix
{\sf diag}$(-1,-1,\dots,-1)$, $|{\sf G}_4|=2$. 
\bea
M_{{\sf G}_4}(t)&=&
\frac{1}{(1-t^2)^n}\sum_{k=0}^{\left[\frac{n}{2}\right]}
\binom{n}{2k}t^{2k}\;. \nonumber \\
P(s,{\sf G}_4) & = &
\left\{ \begin{array}{ll}
\sum_{k=0}^{\left[\frac{n}{2}\right]} \binom{n}{2k}
W(s-2k,{\bf 2}^n) = W(s,{\bf 1}^n),
& s=0 \pmod 2,\\
0, & s \ne 0  \pmod 2.
\end{array}\right.
\label{groupG}
\eea
${\sf G}_4$ is isomorphic as an abstract group to the cyclic group ${\sf
Z}_2$ acting on Euclidean vector space ${\mathbb R}^n$.

It is easy to see that both groups ${\sf G}_2$ and ${\sf G}_4$
acting on ${\mathbb R}^2$ give rise to the same Molien function and
corresponding number of invariants
\bea
M_{{\sf Z}_2}(t)=\frac{1+t^2}{(1-t^2)^2}\;,\;\;\;
P(s,{\sf Z}_2)=\left\{ \begin{array}{ll}
W(s,{\bf 1}^2),
& s=0 \pmod 2,\\
0, & s \ne 0  \pmod 2.
\end{array}\right.
\label{groupPR}
\eea

5. Group ${\sf Q}_{4n}$ generated by the matrices {\footnotesize 
$\left(\begin{array}{cc}
\rho_{2n} & 0 \\
0 & \rho_{2n}^{-1}\end{array}\right)$} and {\footnotesize
$\left(\begin{array}{cc}
0 & i \\
i & 0\end{array}\right)$}, $|{\sf Q}_{4n}|=4n$. 
\bea
M_{{\sf Q}_{4n}}(t) &= &
\frac{1+t^{2n+2}}{(1-t^4)(1-t^{2n})}\;, 
\label{groupQ4n} \\
P(s,{\sf Q}_{4n})&=&\left\{ \begin{array}{ll}
W(\frac{s}{2},\{2,n\}) + W(\frac{s}{2}-n-1,\{2,n\}),
& s=0 \pmod 2,\\
0, & s \ne 0  \pmod 2. \nonumber
\end{array}\right.
\eea
In the case of quaternion group ${\sf Q}_8$ formula (\ref{groupQ4n}) is 
reduced to
\bea
M_{{\sf Q}_8}(t) &= & \frac{1+t^6}{(1-t^4)^2}\;,\;\;\;
P(s,{\sf Q}_8) =\left\{ \begin{array}{ll}
W(s,{\bf 1}^2)/2,
& s=0 \pmod 4,\\
0, & s \ne 0  \pmod 4.
\end{array}\right.
\label{groupQ8}
\eea
More sophisticated examples of the finite groups one can find in Appendices A, B
of the book \cite{Benson93}.
%%%%%%%%%%%%%%%%%%%%%%%%%%%%%%%%%%%%%%%%%%%%%%%%%%%%%%%%%%%%%%%%%%%%%
\section{Conclusion}

1. The explicit expression for restricted partition function $W(s,{\bf 
d}^m)$ and its quasiperiodic components $W_j(s,{\bf d}^m)$ ({\em Sylvester 
waves}) for a set of positive integers ${\bf d}^m = \{d_1, d_2, \ldots, 
d_m\}$ is derived. The formulas are represented as a finite sum over
Bernoulli and Euler polynomials of higher order with periodic 
coefficients.

\noindent
2. Every Sylvester wave $W_j(s,{\bf d}^m)$ satisfies the same 
recursive relation as the whole partition function $W(s,{\bf d}^m)$.

\noindent
3. The application of restricted partition function to the problem of
counting all algebraically independent invariants of the degree $s$
which arise due to action of finite group $G$ on the vector space $V$ over 
the field of complex numbers is discussed.

%\newpage

\section*{Appendices}
\appendix
\renewcommand{\theequation}{\thesection\arabic{equation}}
\section{Symbolic Notation \label{appendix1}}
\setcounter{equation}{0}

The symbolic technique for manipulating sums with binomial coefficients by
expanding
polynomials and then replacing powers by subscripts was developed in
nineteenth century by Blissard.
It has been known as symbolic notation and the classical umbral
calculus \cite{Roman1978}. This notation can be used \cite{Gessel} to
prove
interesting formulas not easily proved by other methods.
An example of this notation is also found in \cite{bat53} in
section devoted to the Bernoulli polynomials $B_k(x)$.

The well-known formulas
$$
B_n(x+y) = \sum_{k=0}^{n} {n\choose k} B_k(x) y^{n-k}, \ \
B_n(x) = \sum_{k=0}^{n} {n\choose k} B_k x^{n-k},
$$
are written
symbolically as
$$
B_n(x+y) = (B(x)+y)^n, \ \ B_n(x) = (B+x)^n.
$$
After the expansion the exponents of $B(x)$ and $B$ are converted into the
orders of the Bernoulli polynomial and the Bernoulli number, respectively:
\be
[B(x)]^k \Rightarrow B_k(x), \ \ \ B^k \Rightarrow B_k.
\label{conv_rule}
\ee
We use this notation in its extended version suggested in
\cite{NorlundMemo} in order to make derivation more clear and
intelligible.
N\"orlund introduced the Bernoulli polynomials of higher order defined 
through
the recursion
\be
B_{n}^{(m)}(x|{\bf d}^m) =
\sum_{k=0}^n  \binom{n}{k} d^k B_k(0) B_{n-k}^{(m-1)}(x|{\bf d}^{m-1}),
\label{Bern_poly_HO_def}
\ee
starting from $B_{n}^{(1)}(x|d_1) = d_1^n B_n(\frac{x}{d_1})$.
In symbolic notation it takes form
$$
B_{n}^{(m)}(x) =  
\left(
d_m B(0) + B^{(m-1)}(x)
\right)^n,
$$
and recursively reduces to more symmetric form
\be
B_{n}^{(m)}(x|{\bf d}^m) =
\left(
x + d_1 \;{}^1\! B(0) +
d_2 \;{}^2\! B(0) + \ldots + d_m \;{}^m\! B(0)
\right)^n =
\left(
x + \sum_{i=1}^m d_i \;{}^i\! B(0)
\right)^n,
\label{Bern_poly_HO_symm}
\ee
where each  $[{}^i \! B(0)]^k$ is converted into $B_k(0)$.
%%%%%%%%%%%%%%%%%%%%%%%%%%%%%%%%%%%%%%%%%%%%%%%%%%%%%%%%%%%%%%%%%%%%%
%\renewcommand{\theequation}{\thesection\arabic{equation}}
\label{appendix2}
\section{Partition function for two primes}
\setcounter{equation}{0}
%We derive the expression for Sylvester waves given in (\ref{sylvp1p2}).
The polynomial part is computed according to (\ref{W_1})
\be
W_1(ap_1p_2,\{p_1,p_2\}) =  
\frac{1}{p_1p_2}
B_{1}^{(2)}(ap_1p_2 + p_1+p_2 |
\{p_1,p_2\})=a+\frac{1}{2}\left(\frac{1}{p_1}+
\frac{1}{p_2}\right)\;.
\label{w1p1p2}
\ee
Two other waves read
\be
W_{p_1}(ap_1p_2,\{p_1,p_2\})=\frac{1}{p_1} \sum_{r=1}^{p_1-1}
\frac{1}{1-\rho_{p_1}^{r}}\;,\;\;\;
W_{p_2}(ap_1p_2,\{p_1,p_2\})=\frac{1}{p_2} \sum_{r=1}^{p_2-1}
\frac{1}{1-\rho_{p_2}^{r}}\;.
\label{w12}
\ee
where we use trivial identity $\rho_{p_1}^{ap_1p_2}=\rho_{p_2}^{ap_1p_2}=1$.
Computation of the sums in (\ref{w12}) 
we start with the identity (see \cite{Vandiver1942})
\be
\prod_{r=0}^{m-1} (x-\rho_m^r) = x^m-1,
\label{identity1}
\ee
and differentiation it with respect to $x$, and division by $x^m-1$
\be
\sum_{r=0}^{m-1} \frac{1}{x-\rho_m^r} = \frac{mx^{m-1}}{x^m-1}.
\label{identity1diff1}
\ee
Subtracting $1/(x-1)$ from both sides of (\ref{identity1diff1}) and taking a
limit at $x \rightarrow 1$ we obtain
\be
\sum_{r=1}^{m-1} \frac{1}{1-\rho_m^r} = \frac{m-1}{2}.
\label{form2}
\ee
Using this result we have for the periodic waves in (\ref{w12})
\be
W_{p_1}(ap_1p_2,\{p_1,p_2\})=\frac{p_1-1}{2p_1}\;,\;\;\;
W_{p_2}(ap_1p_2,\{p_1,p_2\})=\frac{p_2-1}{2p_2}\;.
\label{w12f}
\ee

%\newpage
%%%%%%%%%%%%%%%%%%%%%%%%%%%%%%%%%%%%%%%%%%%%%%%%%%%%%%%%%%%%%%%%%%%%%
\section*{Acknowledgment}
We thank I. M. Gessel for information about Ref. \cite{Gessel}.
The research was supported in part (LGF) by the Gileadi Fellowship
program of the Ministry of Absorption of the State of Israel.

%%%%%%%%%%%%%%%%%%%%%%%%%%%%%%%%%%%%
%%%%%     BIBLIOGRAPHY
%%%%%%%%%%%%%%%%%%%%%%%%%%%%%%%%%%%%

%%%%%%%%%%%%%%%%%%%%%%%%%%%%%%%%%%%%%%%%%%%%%%%%%%%%%%%%%%%%%%%%%%%%%%%%%%%
%%%%%%%%%%%%%%%%%%%%%%%%%%%%%%%%%%%%%%%%%%%%%%%%%%%%%%%%%%%%%%%%%%%%%%%%%%%
\end{document}